\documentclass[12pt,a4paper,twoside]{article}
\usepackage{amsmath,amsfonts,amssymb,a4,enumerate,epsfig,array,theorem}

\newtheorem{theo}{Theorem}
\newtheorem{fact}{Fact}
\newtheorem{prop}{Proposition}[section]
\newtheorem{coro}[prop]{Corollary}
\newtheorem{defin}[prop]{Definition}
\newtheorem{lemma}[prop]{Lemma}

\newcommand{\lebesgue}{\mu}
\newcommand{\tol}{\textbf{tol}}
\newcommand{\id}{\textbf{id}}
\newcommand{\incl}{\textbf{incl}}
\newcommand{\interior}{\operatorname{int}}

\newcommand{\Hsin}{{H_{\sin}}}
\newcommand{\NN}{{\mathbb{N}}}
\newcommand{\ZZ}{{\mathbb{Z}}}
\newcommand{\RR}{{\mathbb{R}}}

\newcommand{\Ss}{{\mathbb{S}}}
\newcommand{\BB}{{\mathbb{B}}}

\newcommand{\HH}{{\mathbb{H}}}

\newcommand{\Nast}{\NN^\ast}
\newcommand{\menos}{\setminus}
\newcommand{\pV}{\partial V}
\newcommand{\pW}{\partial W}

\newcommand{\nobf}{\noindent\bf}

\def\qed{\unskip\nobreak\hfil\penalty50\hskip1.75em\null\nobreak\hfil
$\blacksquare$ {\parfillskip=0pt \finalhyphendemerits=0 \par}\goodbreak}

\begin{document}
\title{Results on infinite dimensional topology 
and applications to the structure of the critical set of
nonlinear Sturm-Liouville operators}
\author{Dan Burghelea, Nicolau C. Saldanha and Carlos Tomei}
\maketitle

\begin{abstract}
We consider  the nonlinear Sturm-Liouville differential operator
$F(u) = -u'' + f(u)$ for $u \in H^2_D([0, \pi])$,
a Sobolev space of functions satisfying Dirichlet boundary conditions.
For a generic nonlinearity $f: \RR \to \RR$ we show that there is
a diffeomorphism in the domain of $F$ converting the critical set $C$ of $F$
into a union of isolated parallel hyperplanes.
For the proof, we show that  the homotopy groups of connected components of
$C$ are trivial
and prove results which permit to replace homotopy equivalences of systems
of infinite dimensional Hilbert manifolds by diffeomorphisms.
\end{abstract}

\section{Introduction}
\label{section:intro}

\footnote{2000 {\em Mathematics Subject Classification}.
Primary 34L30, 58B05, Secondary 34B15, 46T05.
{\em Keywords and phrases} Sturm-Liouville,
nonlinear differential operators,
infinite dimensional manifolds.}
\footnote{D. Burghelea was supported by an NSF grant.
N. Saldanha and C. Tomei acknowledge the support of
CNPq, Faperj and Finep.}
Consider the nonlinear Sturm-Liouville problem
\[ - u''(t) + f(u(t)) = g(t), \qquad u(0) = u(\pi) = 0 \]
and for any smooth nonlinearity $f: \RR \to \RR$
denote by $F$ the differential operator
\begin{align} F: H^2_D([0,\pi]) &\to L^2([0,\pi]). \notag\\
u &\mapsto -u'' + f(u) \notag \end{align}
Here $H^2_D([0,\pi])$ is the Sobolev space of functions $u(t)$ with
square integrable second derivatives which satisfy Dirichlet
boundary conditions $u(0) = u(\pi) = 0$
and $L^2([0,\pi])$ is the usual Hilbert space of square integrable 
functions in $[0,\pi]$.
We are interested in the critical set $C$ of $F$,
\[ C = \{ u \in H^2_D([0,\pi]) \;|\;
DF(u) \textrm{ has 0 as an eigenvalue} \}. \]
Here $DF(u): H^2_D([0,\pi]) \subset L^2([0,\pi]) \to L^2([0,\pi])$
is the Fredholm linear operator $DF(u) v = -v'' + f'(u) v$ of index 0.
Let $\Sigma = \{ m \in \ZZ \;|\; m > 0, -m^2 \in f'(\RR) \}$.
Our main result about $C$ is the following.

\begin{theo}
\label{theo:INTRO}
For tame nonlinearities $f$ (to be defined below),
$C$ is the disjoint union of connected components $C_m$, $m \in \Sigma$,
where each $C_m$ is a smooth contractible hypersurface in $H^2_D([0,\pi])$.
Furthermore, there is a diffeomorphism of $H^2_D([0,\pi])$ to
itself taking $C$ to a union of parallel hyperplanes.
\end{theo}

In the proof of this theorem we shall make use of two results
concerning the topology of infinite dimensional manifolds,
Theorem \ref{theo:A} (proved in Section \ref{section:homoeq})
and Theorem \ref{theo:C} (in Section \ref{section:global}).
The topological theorems are described in a generality which is 
greater than needed in this paper, since we believe that the results are
of independent interest.

\begin{theo}\label{theo:A}
Let $X$ and $Y$ be separable Banach spaces.
Suppose $i:Y \to X$ is a bounded, injective linear map with dense image and 
$M \subset X$ a smooth, closed submanifold of finite codimension.
Then $N = i^{-1}(M)$ is a smooth closed submanifold of $Y$, and 
the restrictions $i:Y \menos N \to X \menos M$ and $i:(Y, N) \to (X,M)$ 
are  homotopy equivalences.
\end{theo}

Our second topological result concerns $H$-manifolds, i.e.,
manifolds modeled on the separable infinite dimensional Hilbert space $H$.

\begin{theo}\label{theo:C}
Suppose $f:(V_1, \pV_1) \to (V_2, \pV_2)$ is a smooth homotopy equivalence
of $H$-manifolds with boundary, $K_2 \subset V_2 \menos \pV_2$ a closed
submanifold of finite codimension and $K_1 = f^{-1}(K_2)$. Suppose also that
$f$ is transversal to $K_2$ and the maps $f:K_1 \to K_2$
and $f: V_1 \menos K_1 \to V_2 \menos K_2$ are homotopy equivalences. 
Then there exists a diffeomorphism $h:(V_1; \pV_1, K_1) \to (V_2; \pV_2, K_2)$,
which is homotopic to $f$ as maps of triples.
\end{theo}

We are able to extend this theorem to a class of Banach spaces
which unfortunately does not contain $C^r_D([0,\pi])$.
We thus cannot obtain a version of Theorem \ref{theo:INTRO}
for such domains of $F$.




We now sketch the main steps in the proof of Theorem \ref{theo:INTRO}.
Let $f : \RR \to \RR$ be a smooth function:
we call $f$ \textit{appropriate} if \textit{either} $f''(0) \ne 0$
\textit{or} the two conditions below hold:
\begin{enumerate}[(a)]
\item{the roots of $f''$ are isolated,}
\item{$f'(0)$ is not of the form $- m^2$, $m \in \ZZ$.}
\end{enumerate}
We call an appropriate function $f: \RR \to \RR$ {\sl tame} iff
$f''(x) \ne 0$ whenever $f'(x)$ is of the form $-m^2$ for some integer $m$.
Notice that tame functions are generic.
In section \ref{section:construct} we introduce
a smooth functional $\omega: H^2_D \to \RR$;
for tame $f$, $\Sigma$ consists of regular values of $\omega$
and $C$ is the union of the non-empty smooth $H$-manifolds
$C_m = \omega^{-1}(\{m\})$.
Theorem \ref{theo:INTRO} follows from Theorem \ref{theo:C}
once we prove that each $C_m$ is contractible.
As any $H$-manifold is an ANR, this fact is asserted by
the following technical result,
which is the core of section \ref{section:construct}.

\begin{prop}
\label{prop:D7}
If $f$ is tame and $C_m$ is non-empty,
then $C_m$ is path-connected and its homotopy groups $\pi_k(C_m)$
are all trivial.
\end{prop}

To prove this proposition,
we observe that the functional $\omega$ smoothly extends to
$\tilde\omega: C^0_D([0,\pi]) \to \RR$ (Lemma \ref{lemma:D1}).
Theorem \ref{theo:A} then shows that the inclusion
$\iota: H^2_D \to C^0_D$ induces a
homotopy equivalence between the levels $C_m \subset H^2_D$ and the levels 
$\tilde C_m = \tilde\omega^{-1}(m\pi) \subset C^0_D, m \in \Sigma$, and we
are left with showing that the manifolds $\tilde C_m$ are
path connected and have trivial homotopy groups:
this is simpler than the similar task for $C_m \subset H^2_D$, 
since we only have to control the continuity of the homotopy
of spheres to a point with respect to the weaker $C^0$ norm. 

\bigskip

We refer to the geometric and topological study
of the set of solutions of $F(u) = g$ (for varying $g$)
as the \textit{geometric approach}.
A pioneering example of the geometric approach applied to PDEs
is the work of Ambrosetti and Prodi on the Laplacian
on a bounded open set $\Omega \subset \RR^n$
with Dirichlet conditions (\cite{AP}),
\[ F_{AP}(u) = - \Delta u + f(u), u|_{\partial \Omega} = 0. \]
In the Ambrosetti-Prodi scenario, the hypotheses on the nonlinearity
are such that the critical set is diffeomorphic to a hyperplane.
Subsequent work (\cite{BP}) then established that $F_{AP}$ is a global fold.
Theorem \ref{theo:INTRO} above is the $n = 1$ case of Ambrosetti-Prodi
but now we consider more general nonlinearities $f$.
For convex $f$, Theorem \ref{theo:INTRO} was proved (\cite{Ruf}, \cite{BT})
by showing that each connected component of $C$ is the graph of a
continuous function from $\Hsin$ to $\RR$, where $\Hsin$ is
the hyperplane of functions in $H^2_D([0,\pi])$ orthogonal to $\sin(t)$.
This result, applied to the nonlinearity $f(u) = u^2/2$,
yields the following rather standard (but not trivial) fact 
in spectral theory:
the set of potentials $u \in H^2_D([0,\pi])$ for which the operator
\[ v \in H^2([0,\pi]) \mapsto -v'' +uv \in L^2([0,\pi]) \]
has $0$ as its $n$-th eigenvalue is a topological hyperplane.
For more general nonlinearities, however,
we were not able to make spectral theory work for us.
Our hypotheses do not demand that the nonlinearity $f$
have a prescribed asymptotic behavior at infinity.

A more elementary version of this approach has been exploited
in \cite{MST1} to show that the critical set of the operator
\begin{align}
                       F_1: H^1(\Ss^1)  &\to L^2(\Ss^1)       \notag\\
                             u &\mapsto   -u' + f(u)         \notag
\end{align}
is either empty or a hyperplane. In this case, the critical
set is the zero level of a Nemytskii operator,
\begin{align}
                        \varphi: H^1(\Ss^1) &\to \RR            \notag\\
                                  u &\mapsto \int_{\Ss^1} f'(u) \notag
\end{align}
whose contractibility was shown (\cite{MST2})
by means of ergodic-like arguments,
robust enough to admit extensions to functionals
from spaces of functions acting in domains 
in higher dimensions and taking values on $\RR^n$. 
To show contractibility of the connected components $C_m$ of $F$,
however, we recur constantly to Sturm oscillation,
and this is the main reason why our proof does not appear to extend
to operators on functions in many variables.
Notice that Lemma 5.3 in \cite{MST1} is a corollary of our Theorem \ref{theo:C}.

In order to provide global geometric information about the operator $F_1$,
the authors of \cite{MST1} considered
the stratification of the critical set into 
Morin singularities of different types. Generically, the
singularities of $F$ are also of Morin type, but we do not explore the
matter further in this paper.

\section{Homotopy equivalence}
\label{section:homoeq}

The aim of this section is to prove Theorem \ref{theo:A}.
As in the statement of the theorem, $X$ and $Y$ are separable Banach spaces,
$i: Y \to X$ is an injective bounded linear map with dense image
and $M \subset X$ is a closed smooth submanifold of finite codimension $k$.
Let $P$ be a compact smooth manifold with boundary of dimension $r + k$;
$P$ has a fixed but arbitrary Riemannian metric.

\begin{defin}
The smooth function $f: P \to X$
is called a {\em smooth $M$-proper embedding} if
$f$ is a smooth embedding,
$f(\partial P) \cap M = \emptyset$ and
$f$ is transversal to $M$.
\end{defin}

The set $Q = f^{-1}(M)$ is a smooth compact manifold (with no boundary)
of dimension $r$.

\begin{lemma}
\label{lemma:C1stable}
If $f: P \to X$ is a smooth $M$-proper embedding then there exists
$\epsilon > 0$ such that any smooth map $g: P \to X$
with $\|f - g\|_{C^1} < \epsilon$ is a smooth $M$-proper embedding.
For this $\epsilon$,
if $\|f - g\|_{C^1} < \epsilon$ then $Q_g = g^{-1}(M)$
is diffeomorphic to $Q$ and there exists a smooth embedding $\theta: Q \to P$
homotopic to the inclusion $Q \subset P$ and
such $\theta(Q) = Q_g$.
Moreover, if $S \subset P$ is a compact submanifold and $f|_S = g|_S$
then $\theta$ can be chosen to be the identity on $S \cap Q$ and
the homotopy to be relative to $S \cap Q$.
\end{lemma}

Notice that the lemma also holds if $Q = \emptyset$.

{\nobf Proof:}
The existence of $\epsilon$ follows from the fact that all of the following
properties are open in $f$ in the $C^1$ topology:
being an embedding, $f(\partial P) \cap M = \emptyset$
and $f$ is transversal to $M$.
Let $g_t: P \to X$ be defined by $g_t(p) = (1-t)f(p) + tg(p)$;
clearly $g_0 = f$ and $g_1 = g$ and by the previous remark
$g_t$ is a smooth $M$-proper embedding for all $t \in [0,1]$.
Let $G: P \times [0,1] \to X \times [0,1]$ be $G(p,t) = (g_t(p),t)$
and let $\tilde Q = G^{-1}(M \times [0,1])$:
$\tilde Q$ is a compact manifold with two boundary components,
$Q \times \{0\}$ and $Q_g \times \{1\}$.
The function $\pi: \tilde Q \to [0,1]$ defined by $\pi(p,t) = t$
is a submersion.
Notice that $(S \cap Q) \times [0,1] \subseteq \tilde Q$.
Construct on $\tilde Q$ a tangent vector field $\alpha$
such that $\alpha(p,t) = (0,1)$ for $p \in S \cap Q$
and $D\pi(p,t) \cdot \alpha(p,t) = 1$ for all $(p,t) \in \tilde Q$
(it is easy to construct such a vector field $\alpha$ in a neighborhood
of a point $(p,t)$; use partitions of unity to define it on all $\tilde Q$).
Integrating this vector field yields $\theta$ and the desired homotopy.
\qed

Let $X$, $Y$ and $i: Y \to X$ be as above and let $P$
be a compact manifold with boundary.
Let $C^1(P,X)$ (resp. $C^1(P,Y)$) be the metric space of $C^1$ functions
from $P$ to $X$ (resp. $Y$) with the $C^1$ metric.
Similarly, let $C^1_c(\RR^k,X)$ (resp. $C^1_c(\RR^k,Y)$)
be the normed vector spaces of $C^1$ functions from $\RR^k$ to $X$ (resp. $Y$)
with compact support with the $C^1$ norm.
Define $i^\ast: C^1(P,Y) \to C^1(P,X)$
and $i^\ast: C^1_c(\RR^k,Y) \to C^1_c(\RR^k,X)$ by composition.

\begin{lemma}
\label{lemma:C1dense}
The images of $i^\ast: C^1(P,Y) \to C^1(P,X)$
and $i^\ast: C^1_c(\RR^k,Y) \to C^1_c(\RR^k,X)$
are dense in $C^1(P,X)$ and $C^1_c(\RR^k,X)$, respectively.
\end{lemma}

{\nobf Proof:}
We first prove the lemma for $\RR^k$ by induction on $k$
(the case $k=0$ is trivial).
Let $f: \RR^{k+1} \to X$ be a $C^1$ function with compact support
and $\epsilon > 0$ a real number.
We may assume without loss of generality that $f$ is smooth
(take a convolution with a smooth bump) and that the support of $f$
is contained in $(0,1)^{k+1}$.
We want to construct $\tilde f: \RR^{k+1} \to Y$
such that $d_{C^1}(f,i \circ \tilde f) < \epsilon$
and such that the support of $\tilde f$ is also contained in $(0,1)^{k+1}$.
Take $\delta = 1/N > 0$ such that if $v, v' \in \RR^{k+1}$, $d(v,v') < \delta$,
then $f$ and all partial derivatives of $f$ of order 1 or 2 differ by
at most $\epsilon/16$ between $v$ and $v'$.
Assume furthermore that the support of $f$
is contained in $(\delta, 1 - \delta)^{k+1}$.
Let $g = \partial f/\partial x_{k+1}$ and
consider the functions $g_0, g_1, \ldots, g_N: \RR^k \to X$ defined by
\[ g_j(x_1,x_2,\ldots,x_k) = g(x_1,x_2,\ldots,x_k,j/N); \]
notice that $g_0 = g_1 = g_{N-1} = g_N = 0$.
By induction hypothesis,
we may pick $\tilde g_0, \ldots, \tilde g_N: \RR^k \to Y$
with $d_{C^1}(g_j, i \circ \tilde g_j) < \epsilon/16$ and with supports
contained in $(0,1)^k$;
take $\tilde g_0 = \tilde g_1 = \tilde g_{N-1} = \tilde g_N = 0$.
Now define $\tilde g: \RR^{k+1} \to Y$ by
\[ \tilde g(x_1,x_2,\ldots,x_k,j/N) = \tilde g_j(x_1,x_2,\ldots,x_k) \]
and by affine interpolation for other values of $x_{k+1}$.
Clearly, the distances
$d_{C^0}(g,i\circ\tilde g)$,
$d_{C^0}(\partial g / \partial x_1, i \circ \partial \tilde g / \partial x_1)$,
\dots,
$d_{C^0}(\partial g / \partial x_k, i \circ \partial \tilde g / \partial x_k)$
are all smaller than $\epsilon/4$.
Therefore, the function $\tilde h: \RR^n \to Y$ defined by
$\tilde h(x_1,\ldots,x_n) = \int_0^1 \tilde g(x_1,\ldots,x_k,t) dt$
satisfies $d_{C^1}(h,0)  < \epsilon/4$.
Let $\phi: \RR \to \RR$ be a smooth non-negative function
with support contained in $(0,1)$,
integral equal to 1 and $d_{C^1}(\phi,0) < 3$.
Then the function
\[ \tilde f(x_1,\ldots,x_k,x_{k+1}) =
\int_0^{x_{k+1}} (\tilde g(x_1,\ldots,x_k,t) -
\phi(t) \tilde h(x_1,\ldots,x_k)) dt \]
satisfies all the requirements.

We now prove the lemma for a compact manifold $P$.
Take a finite open cover of $P$ by disks
and a corresponding smooth partition of unity.
In order to approximate $f$ it suffices to approximate each product of $f$
by a function in the partition of unity,
provided the support of the approximation
is still contained in the corresponding open set.
But that is precisely what we did in the previous case.
\qed

Let $M$ be a submanifold of finite codimension $k$
of a separable Banach space $X$.
A closed tubular neighborhood of $M$ is a 0-codimensional
smooth embedding $\phi: D(\xi) \to X$, where $D(\xi)$ is the closed unit
disk bundle of a smooth $\RR^k$-bundle over $M$ so that
\begin{enumerate}[(1)]
\item{$\phi$ restricted to the 0-section is the inclusion
$M \hookrightarrow X$,}
\item{$\phi(D^0(\xi))$ is an open subset of $X$,}
\item{$\phi(D(\xi))$ is a closed subset of $X$.}
\end{enumerate}
Here, $D^0(\xi)$ is the open unit disk bundle. Clearly,
$\phi(\partial D(\xi))$ is a codimension 1 smooth submanifold of $X$.
It is a well known fact (\cite{Lang} or \cite{Elw})
that finite codimensional Banach submanifolds
of a separable Banach space admit closed tubular neighborhoods.

{\nobf Proof of Theorem \ref{theo:A}:}
It suffices to prove the two following facts.
Let $S$ and $Q$ be compact manifolds with no boundary.
If $i_0: S \to Q$, $i_1: S \to Y \menos N$ (resp. $i_1: S \to N$)
and $i_2: Q \to X \menos M$ (resp. $i_2: Q \to M$)
are embeddings with $i \circ i_1 = i_2 \circ i_0$
then there exists $u: Q \to Y \menos N$ (resp. $u: Q \to N$)
such that $u \circ i_0$ is homotopic to $i_1$ and
$i_2$ is homotopic to $i \circ u$.

In the first case we use Lemma \ref{lemma:C1stable} and
Lemma \ref{lemma:C1dense} to obtain $\tilde i_2: P \to X \menos M$
near $i_2$ and of the from $\tilde i_2 = i \circ u$,
for $u: P \to Y \menos N$, proving the first claim.

Consider a closed tubular neighborhood of $M$ in $X$.
If $i_2: Q \to M$ is an embedding then the pull-back of the tubular
neighborhood is a bundle $P \to Q$.
We have a smooth $M$-proper embedding also called $i_2: P \to X$
where $P$ is a compact manifold with boundary;
the dimension of $P$ is $r+k$ where $k$ is the codimension of $M$
and $r$ is the dimension of $Q$.
Use this construction in the second case to define $P$,
again use both lemmas to obtain $\tilde i_2: P \to X$ near $i_2$,
a smooth $M$-proper embedding of the form $\tilde i_2 = i \circ u_1$,
where $u_1: P \to Y$ a smooth $N$-proper embedding.
The homotopy and the function $\theta$ constructed in
Lemma \ref{lemma:C1stable}
now obtain $u: Q \to Y$ and the desired homotopies.
\qed

\section{Global changes of variable}
\label{section:global}


In this section we prove Theorem \ref{theo:C}.
We will write $\HH$ for the infinite dimensional separable real Hilbert space
and we call a $\HH$-manifold (resp. $\HH$-manifold with boundary)
a Hausdorff paracompact smooth manifold with local model $\HH$
(resp. $\HH \times [0,+\infty)$).

\begin{prop}\label{prop:B}
Suppose $f: (V, \partial V) \to (W, \partial W)$ is a homotopy equivalence
between two smooth $\HH$-manifolds with boundary

\begin{enumerate}
\item{There exists a diffeomorphism $h:(V, \partial V) \to (W, \partial W)$
so that $h$ and $f$ are homotopic maps of pairs.}
\item{If $h^{\partial}: \partial V \to \partial W$ is a diffeomorphism 
homotopic (resp. equal) to $f^{\partial} = f|_{\pV} : \pV \to \pW$,
then one can extend $h^{\partial}$ to a diffeomorphism
$h: (V, \pV) \to (W, \pW)$ homotopic (resp. relative homotopic) to $f$.}
\end{enumerate}
\end{prop}

Here, an $\HH$-manifold is modeled on the
separable infinite dimensional real Hilbert space.
The proof of Proposition \ref{prop:B} is based on
the following known results on Hilbert manifolds.

\begin{fact}\label{fact:2.1} \textrm{\cite{BK}, \cite{Burg}} 
If $f: M \to N$ is a homotopy equivalence between two $\HH$-manifolds, there
exists $h: M \to N$, a diffeomorphism which is homotopic to $f$.
\end{fact}

An isotopy between diffeomorphisms $h_0, h_1 : M \to N$ 
is a diffeomorphism $h:\RR \times M \to \RR \times N$ taking $(t,m)$
to a point of the form $h(t,m) = (t,h_t(m))$ 
so that $h_t = h_0$ if $t \le 0$ and $h_t = h_1$ if $t \ge 1$.

\begin{fact}\label{fact:2.2} \textrm{\cite{Burg}}
Homotopic diffeomorphisms $h_0, h_1: M\to N$  between $\HH$-manifolds are
isotopic.
\end{fact}

\begin{fact}\label{fact:2.3} \textrm{\cite{BK}, \cite{Elw}}
Given two homotopic closed embeddings of infinite codimension $\ell_i:V \to M$,
$i=0,1$, with $V$ and $M$ $\HH$-manifolds, there exists an isotopy
of diffeomorphisms $h: \RR \times M \to  \RR \times M$ such that $h_t$ is the
identity for $t \le 0$, $h_t$ is constant for $t \ge 1$ and
$h_1 \circ \ell_0 = \ell_1$.
\end{fact}

Let $M$ be an $\HH$-manifold and $V \subset M$ be a closed 
$\HH$-submanifold. As in Section \ref{section:homoeq}
define a closed tubular neighborhood of $V$ to be a 0-codimensional 
smooth embedding $\phi: D(\xi) \to M$, where $D(\xi)$ is the closed unit 
disk bundle of a smooth vector bundle with the same properties (1),
(2) and (3) as before.
Now, the fiber of $\xi$ is isomorphic to $\RR^k$
if the codimension of $V$ in $M$ is $k$
and $\HH$ if the codimension is infinite.
If the fiber is isomorphic to $\HH$ then
$D(\xi)$ is necessarily a trivial bundle since
the linear group of invertible bounded operators on $\HH$ is contractible
(\cite{Ku}).

\begin{fact}\label{fact:2.4}\textrm{\cite{BK}}
Let $V \subset M $ as above and $\phi_i: D(\xi) \to M$, $i=0,1$, be two
closed tubular neighborhoods. Then there exists an isotopy of diffeomorphisms 
$h: \RR \times M \to \RR \times M$ such that $h_t$ is the identity
for $t \le 0$, $h_t$ is constant for $t \ge 1$,
$h_1 \circ \phi_0 = \phi_1 \circ \theta$, $\theta: D(\xi) \to D(\xi)$
a (vector) bundle isomorphism
which can be taken to be identity when the fiber of $\xi$ is $\HH$
and $h_1|_{V}$ is the identity.
\end{fact}

The hypothesis of finite codimension implies ${K_1}$ and ${K_2}$ of infinite
dimension. A similar conclusion remains true for ${K_1}$ and ${K_2}$ of 
infinite dimension and infinite codimension,
and the proof in this case is a straightforward consequence
of results in \cite{BK} and \cite{Burg}.

{\nobf Proof of Proposition \ref{prop:B}}
Let $V_0 = V$ and $V_1 = W$.
Set $K_i$, $i=0,1$, to be $\HH$-manifolds diffeomorphic
to $\pV_i$ via diffeomorphisms $\psi_i: \pV_i \to K_i$.
Let $B_i = K_i \times D^{\infty}$
and $K_i^0 = K_i \times \{0\} \subset B_i$. 
By Fact \ref{fact:2.1}, there exists diffeomorphisms
$\theta_i: \pV_i \to K_i \times \Ss^{\infty}$, such that
$\theta_i$ is homotopic to the map $\tilde \psi_i$ taking $v \in \pV_i$ to 
$(\psi_i(v),x) \in K_i \times \Ss^\infty$, where $x$ is
an arbitrary but fixed element in $\Ss^\infty$.
Also, define the smooth $\HH$-manifolds $\tilde V_i = V_i \cup_{\theta_i} B_i$.

\begin{figure}[ht]
\begin{center}
\epsfig{height=4cm,file=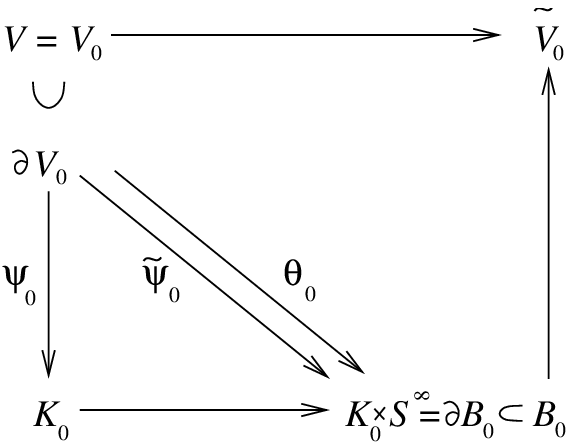} \qquad \epsfig{height=4cm,file=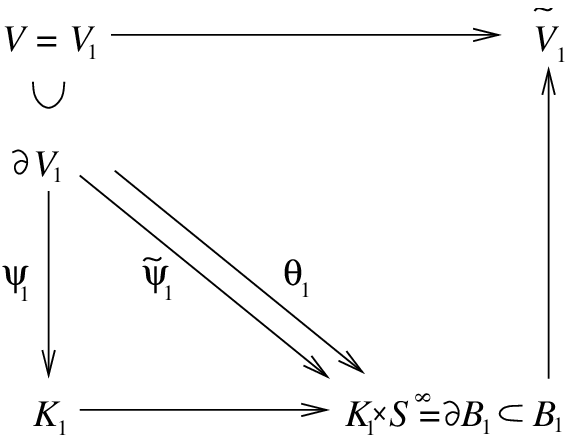}
\end{center}
\label{fig:dia1}
\end{figure}

The inclusions $(\tilde V_i, K_i^0) \hookrightarrow (\tilde V_i, B_i)$ 
and $(V_i, \pV_i) \hookrightarrow (\tilde V_i, B_i)$ are
homotopy equivalences of pairs. Also, by hypothesis, 
$f : (V_0, \pV_0) \to (V_1, \pV_1)$ is a homotopy equivalence of pairs. 
Hence, there exists $\tilde f : (\tilde V_0, K_0^0) \to (\tilde V_1, K_1^0)$
which makes the diagram below homotopy commutative.
From Fact \ref{fact:2.1}, 
let ${\tilde f}': \tilde V_0 \to \tilde V_1$
be a diffeomorphism which is homotopic to $\tilde f$. Notice that
the behavior of ${\tilde f}'$ at $K_0^0$ or $\pV_0$ is not controlled.
Let $\ell : K_0 = K_0^0 \to {K_1} = K_1^0$ be a diffeomorphism
homotopic to $\tilde f|_{K_0^0}$. Notice that $f|_{\pV_0}: \pV_0 \to \pV_1$ is
homotopic to both $\psi_1^{-1} \circ \ell \circ \psi_0$ and to
$\theta_1^{-1} \circ (\ell \times \id) \circ \theta_0$.

\begin{figure}[ht]
\begin{center}
\epsfig{height=5cm,file=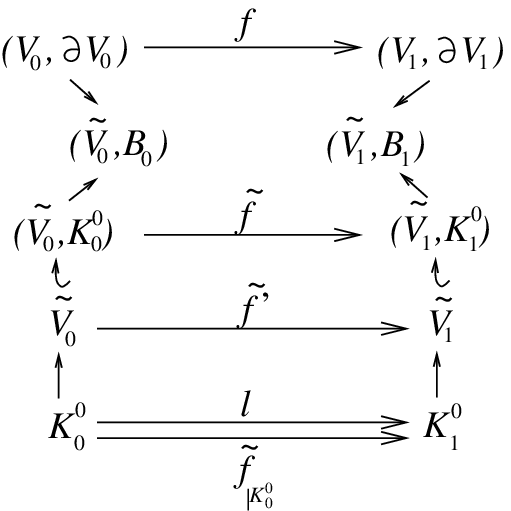}
\end{center}
\label{fig:dia2}
\end{figure}

Let $\ell_0: K_0 \to \tilde V_1$ be the composition of $\ell$ with the 
inclusion ${K_1} = {K_1}^0 \subset \tilde V_1$. 
Also, let $\ell_1: K_0 = K_0^0 \to \tilde V_1$ 
be the restriction of ${\tilde f}'$ to $K_0^0$. Notice that both maps
embed $K_0$ in $\tilde V_1$ as an infinite codimensional submanifold.
By Fact \ref{fact:2.3}, there exists an isotopy 
$h':\RR \times  \tilde V_1 \to \RR \times \tilde V_1$ with 
$h'_t = \id$ for $t \le 0$, which is independent of $t$ for $t \ge 1$ and 
such that $h'_1 \circ \ell_0 = \ell_1$. The diffeomorphism 
${\tilde f}'' = (h'_1)^{-1} \circ {\tilde f}': \tilde V_0 \to \tilde V_1$
is clearly homotopic to $\tilde f$. Also, ${\tilde f}''$ takes
$K_0^0$ to ${K_1}^0$, and ${\tilde f}''$ is then homotopic to $\tilde f$ as a 
map of pairs, from $(\tilde V_0, K_0^0)$ to $(\tilde V_1, {K_1}^0)$.

Let $\phi_0, \phi_1: {K_1}^0 \times D^\infty \to \tilde V_1$ be tubular 
neighborhoods given by $\phi_0 = {\tilde f}'' \circ \incl \circ
(\ell^{-1} \times \id)$, where $\incl: {K_1} \times D^\infty \to \tilde V_1$
is the inclusion, and $\phi_1$ is the inclusion of $B_1$ in $\tilde V_1$.
By Fact \ref{fact:2.4}, there is an isotopy 
$h'': \RR \times \tilde V_1 \to \RR \times \tilde V_1$ with
$h''_t = \id$, for $t \le 0$, which is independent of $t$ for
$t \ge 1$ and such that $h_1'' \circ \phi_0 = \phi_1$.
Finally, let ${\tilde f}''' = h_1'' \circ {\tilde f}''$. Clearly, 
${\tilde f}'''$ restricts to a diffeomorphism of pairs from 
$(V_0, \pV_0)$ to $(V_1, \pV_1)$, homotopic to $f$. This finishes item (a).

Item (b) is a straightforward consequence of item (a), Fact \ref{fact:2.2}
and the existence of collar neighborhoods for the boundary of an $\HH$-manifold.
\qed

Before we proceed with the proof of Theorem \ref{theo:C}
we describe a construction which is sometimes refered to as
\textit{cutting a manifold along a submanifold}.

Let $(V, \partial V)$ be an $\HH$-manifold with boundary. 
A proper finite codimension submanifold is a closed subset $K \subset V$
such that
\begin{enumerate}
\item{$K \subset V \menos \partial V$,}
\item{$K$ is a closed, finite codimensional smooth submanifold of 
$V \menos \partial V$.}
\end{enumerate}
The {\emph{ normal bundle}} $\nu_K$ is the quotient $T(V)|_K / T(K)$.
Denote by $E(\nu_K)$ the total space of $\nu_K$ and by 
$S(\nu_K)$ the corresponding sphere bundle $(E(\nu_K) \menos K) / \RR^+$,
where $\RR^+$ acts on $E(\nu_K) \menos K$ by multiplication. Also, denote
by $D(\nu_K)$ the fiberwise compactification of $E(\nu_K)$ obtained by adding 
a point at infinity for each half-line from the origin.
Clearly, 
$D(\nu_K)$ and $S(\nu_K)$ are isomorphic to the closed unit disk bundle
and to the unit sphere bundle of $\nu_K$
provided $\nu_K$ is equipped with a smooth fiberwise scalar product,
and we may therefore consider $S(\nu_K)$ as the boundary of $D(\nu_K)$. 

Let $\psi_D: D(\nu_K)\to V$ be a closed tubular
neighborhood of K.
Let $\psi_S: S(\nu_K) \times [0,1) \to V$ be defined
by $\psi_S(n , t) = \psi_D(tn)$.
The  {\emph{cutting}} of $(V, \partial V)$ {\emph{along}} $K$ 
is an $\HH$-manifold with boundary (actually a bordism) 
$(\tilde V, \partial_+ \tilde V, \partial_- \tilde V)$, together with
a canonical smooth map 
$p :  (\tilde V, \partial_+ \tilde V, \partial_- \tilde V)
\to (V, \partial V)$.
The manifold $\tilde V$ as a set is the disjoint union of
$V \menos K$ and $S(\nu_K)$; 
more precisely, $\tilde V = (V \menos K) \cup_{\psi_S} (S(\nu_K) \times [0,1))$.
Set $p$ to be the identity on $V \menos K$
and equal to the
bundle projection  $S(\nu_K)\to K$ on $K$.
From Fact \ref{fact:2.4},
the smooth structure thus defined on $\tilde V$ actually does not depend 
on the choice of tubular neighborhood $\theta$. 
Now take $\partial_+ \tilde V = \partial V$
and $\partial_- \tilde V = S(\nu_K)$.

Consider $\check V = \tilde V \cup D(\nu_K)$,
obtained by identifying $\partial_- \tilde V$
and $\partial D(\nu_K)$, both equal to $S(\nu_K)$. The pair 
$(\check V, \partial \check V = \partial_+  \tilde V)$
is a smooth $\HH$-manifold
with boundary and $K \subset D(\nu_K) \subset\check V$
is a proper smooth embedding
of finite codimension.
There exists a well defined class of {\emph{ thickening diffeomorphisms}}
(all isotopic) 
$\check\theta:(\check V, \partial_+ \check V, K) \to (V, \partial V, K)$
so that the restriction of $\check \theta$ to $K$ is the identity. 

The constructions above are functorial in the following sense. 
Let
\[f: (V_1, \partial V_1, {K_1}) \to (V_2, \partial V_2, {K_2})\]
be a smooth map 
so that $f$ is transversal to ${K_2}$ and ${K_1} = f^{-1}({K_2})$. 
Such a map induces a map 
$\tilde f: (\tilde V_1, \partial_+ \tilde V_1, \partial_- \tilde V_1) \to
(\tilde V_2, \partial_+ \tilde V_2, \partial_- \tilde V_2)$ and bundle maps
$D(\nu_{K_1})\to D(\nu_{K_2})$ sending $S(\nu_{K_1})$ into $S(\nu_{K_2})$.
The restriction of $\tilde f$  to $\partial_- \tilde V_1$
is the induced bundle map
from $S(\nu_{K_1})$ to $S(\nu_{K_2})$; on each fiber, this map is projective 
(i.e., $\tilde f$ restricted to $\partial_- \tilde V_1$ comes from 
a vector bundle map $E(\nu_{K_1}) \to E(\nu_{K_2})$).
Such maps $f$ will be called {\emph{ morphisms}}. We consider morphisms
$f: (V_1, \partial V_1, {K_1}) \to (V_2, \partial V_2, {K_2})$
with the additional
property that the maps $f : V_1 \to V_2$, $f: \partial V_1 \to \partial V_2$,
$f : {K_1} \to {K_2}$, $f: V_1 \menos {K_1} \to V_2 \menos {K_2}$
are all homotopy
equivalences. This clearly implies that
$\tilde f: (\tilde V_1, \partial_+ \tilde V_1, \partial_- \tilde V_1) \to
(\tilde V_2, \partial_+ \tilde V_2, \partial_- \tilde V_2)$ is a homotopy
equivalence of triples.

{\nobf Proof of Theorem \ref{theo:C}:}
Start with $f: (V_1, \partial V_1, {K_1}) \to (V_2, \partial V_2, {K_2})$
which induces the homotopy equivalence 
$\tilde f: (\tilde V_1, \partial_+ \tilde V_1, \partial_- \tilde V_1) 
\to (\tilde V_2, \partial_+ \tilde V_2, \partial_- \tilde V_2)$ and
the bundle map $D(f) : D(\nu_{K_1}) \to D(\nu_{K_2})$. 
Notice that the restriction $f^K = f|_{{K_1}}$ is a homotopy equivalence,
and therefore so is $D(f)$. 

By Fact \ref{fact:2.1}, $f^K$ is homotopic to a diffeomorphism
$h^K: {K_1} \to {K_2}$ by a homotopy $f_t^K$, where $f_0^K = f^K$ and
$f_1^K = h^K$. Since $f^K$ comes from a bundle map $D(f)$, the homotopy
$f_t^K$ lifts to a homotopy of bundle maps $D(f_t^K)$, and since 
$D(f_1^K)$ induces on the base a diffeomorphism $h^K$, $D(f_1^K)$
itself is a diffeomorphism. 


This shows that $\tilde f$ is a homotopy equivalence between 
$\partial_- \tilde V_1$ and $\partial_- \tilde V_2$, and therefore
a homotopy equivalence between the boundaries of $\tilde V_1$ and 
$\tilde V_2$. We may then apply Proposition \ref{prop:B} to obtain
a homotopic diffeomorphism $\tilde h$ between the triples. Also,
the restriction of $D(f_1^K)$ to $S(\nu_{K_1}) = \partial_- \tilde V_1$
is a diffeomorphism homotopic to the restriction of $\tilde h$.
By Fact \ref{fact:2.2}, there exists an isotopy between these
two diffeomorphisms, and we may therefore glue them in order to obtain
the desired diffeomorphism.
\qed

\begin{coro}
\label{coro:difeo}
Let $X$ and $Y$ be a separable Hilbert spaces and
$i: Y \to X$ an injective bounded linear map with dense image.
If $M$ is a finite codimensional closed submanifold of $X$
then $N = i^{-1}(M)$ is a finite codimensional closed submanifold of $Y$
and there exists a diffeomorphism $h:(Y,N) \to (X,M)$
homotopic to $i:(Y,N) \to (X,M)$. 
\end{coro}

{\nobf Proof:}
From Theorem \ref{theo:A},
$i: Y \menos N \to X \menos M$ and $i: (Y,N) \to (X,M)$
are homotopy equivalences.
We now apply Theorem \ref{theo:C} with $V_1= Y$ and $V_2 = X$,
$\pV_1 = \pV_2 = \emptyset$, ${K_2} = M$, ${K_1} = N$ and $f = i$.
\qed

\section{Constructing the homotopy}
\label{section:construct}

Let $X = C^0_D([0,\pi])$, $Y = H^2_0([0,\pi])$ and $Z = L^2([0,\pi])$.
Consider the operator
\begin{align}
F: Y &\to Z \notag\\
u &\mapsto -u'' + f(u)\notag
\end{align}
with derivative at $u$ given by
\begin{align}
DF(u): Y &\to Z. \notag\\
w &\mapsto -w'' + f'(u) w\notag
\end{align}
We recall some facts from the theory of second order
differential equations (\cite{CL}).
Given $u \in Y$, let $v(u,\cdot)$ be the solution in $[0,\pi]$ of
\[ -v''(u,t) + f'(u(t)) v(u,t) = 0,
\qquad v(u,0) = 0, \quad v'(u,0) = 1. \]
The characterization of critical points of $F$ follows from
standard Fredholm theory of Sturm-Liouville operators:
$u$ is in $C$ if and only if the kernel of $DF(u)$ is non-trivial.
Thus, a point $u$ belongs to the critical set $C$ of $F$
if and only if $v(u,\pi) = 0$;
in this case, by the simplicity of the spectrum of $DF(u)$,
$\ker DF(u)$ is spanned by $v(u,\cdot)$.

Let $\omega: Y \times [0,\pi] \to \RR$ be the continuously defined
argument of $(v'(u,t), v(u,t))$, with $\omega(u,0) = 0$.
One can show that
$u$ is critical if and only if $\omega(u,\pi) = m \pi$, for $m \in \ZZ$;
moreover, $m$ has to be positive since $\omega(u,t) > 0$ for all $t > 0$.
Let $C_m = \{ u \in Y\;|\;\omega(u,\pi) = m\pi\}$:
the critical set $C$ of $F$
is the (disjoint) union of $C_m$, $m \in \Nast$.
A standard computation shows that, for any given $t$,
$\omega(u,t)$ is smooth as a function of $u \in Y$ and we have
\[ \frac{\partial}{\partial u} \omega(u,t) \cdot \varphi =
- \frac{1}{(v(u,t))^2 + (v'(u,t))^2}
\int_0^t f''(u(s)) \varphi(s) (v(u,s))^2 ds. \]
where $\frac{\partial}{\partial u} \omega(u,t)$ denotes the differential
of $\omega(\cdot,t): Y \to \RR$ and 
$\frac{\partial}{\partial u} \omega(u,t) \cdot \varphi$
the value of this differential on the element $\varphi \in Y$.
Let $f$ be appropriate: 
using the formula above one can see that
$m\pi$ is a regular value for the real valued function $\omega(\cdot, \pi)$,
and therefore
the sets $C_m$ are either empty or smooth manifolds of codimension 1.

When considering a given $C_m$
we use the more convenient {\sl $m$-argument} at $u$: 
$\omega_m: Y \times [0,\pi] \to \RR$ is the argument of $(v', mv)$,
i.e., it satisfies
\[ v'(u,t) \tan(\omega_m(u,t)) = m v(u,t), \quad \omega_m(u,0) = 0. \]
It is easy to see that $m \pi$ is a regular value of $\omega_m$,
as it is of $\omega$.
The advantage of this definition is that if $f'(u(t)) = -m^2$
for $t$ in some interval then $\omega_m(u,\cdot)$ is a linear map
of slope $m$ in this interval.

We will later consider {\it local $m$-arguments}.
More precisely, given $u$ we solve the differential equation
\[ - \hat v'' + f'(u) \hat v = 0,\qquad 
\hat v(t_0) = a_0,\quad \hat v'(t_0) = b_0 \]
and set $\hat\omega_m(t)$ to be the argument of $(\hat v', m\hat v)$.
The notation $\hat\omega_m(u,t)$ leaves unspecified the values of
$t_0, a_0$ and $b_0$:
these values will always be specified in the context.

We list some properties relating $u$ and $\omega_m$,
most of which are simple or standard.
The obvious adaptations to $\hat\omega_m$ will be frequently
left to the reader.

\begin{lemma}
\label{lemma:D1}
The $m$-argument $\omega_m$ satisfies
\begin{equation}
\omega_m'(u,t) =
m - \frac{m^2 + f'(u(t))}{m} \sin^2(\omega_m(u,t)),
\qquad \omega_m(u,0) = 0 \tag{$\ast$}
\end{equation}
(here $\omega_m'(u,t) = \frac{d}{dt} \omega_m(u,t)$);
\begin{enumerate}[(a)]
\item{ $\omega_m'(u,t) = m$ if and only if
either $\omega_m(u,t) = j \pi$, $j \in \ZZ$, or $f'(u(t)) = -m^2$; }
\item{ if $f'(u(t)) < -m^2$ (resp., $f'(u(t)) > -m^2$)
then $\omega_m'(u,t) \ge m$ (resp., $\omega_m'(u,t) \le m$); }
\item{ the differential equation $(\ast)$ defines
$\omega_m(\cdot,\pi)$ as a smooth function on $C^0([0,\pi])$
which is $L^1$-continuous on bounded sets; }
\item{ for any given $t$, the differential equation $(\ast)$ implies that
\[ \frac{\partial}{\partial u} \omega_m(u,t) \cdot \varphi =
- \frac{m}{(m v(u,t))^2 + (v'(u,t))^2}
\int_0^t f''(u(s)) \varphi(s) (v(u,s))^2 ds; \] }
\item{ if $u$ is smooth and not flat at $t_0$,
$f''(u(t_0)) \ne 0$ and $\omega'_m(u,t_0) = m$
then $\omega'_m(u,\cdot)$ is not flat at $t_0$. }
\end{enumerate}
\end{lemma}

Item (c) of this result entitles us to define the extension
$\tilde\omega_m: X \to \RR$
and implies that the maps $\omega$ and $\omega_m$
can be extended to smooth maps on $X$
and one can consider the codimension one smooth manifolds 
$\tilde C_m = \{ u \in X\;|\;\omega(u,\pi) = m\pi\}$.
From Corollary \ref{coro:difeo},
the inclusion $i: (Y,C_m) \to (X,\tilde C_m)$ is homotopic
to a homeomorphism.
From now on we drop the tilde both on $\tilde C_m$ and $\tilde\omega_m$:
thus, our basic functional space will be $X = C^0_D([0,\pi])$.

Recall that a smooth function $\zeta: \RR \to \RR$ is {\sl flat}
at $t_0$ if its Taylor expansion at $t_0$ is constant.

{\nobf Proof of Lemma \ref{lemma:D1}: }
We only prove the last item.
Assume by contradiction that $\omega'_m(u,\cdot)$ is flat at $t_0$.
Since $f''(u(t_0)) \ne 0$, the function $f' \circ u$ is not flat at $t_0$.
In the differential equation
\[ \omega_m'(u,t) - m =
- \frac{m^2 + f'(u(t))}{m} \sin^2(\omega_m(u,t)) \]
the left hand side is zero and flat at $t_0$ and
the right hand side is a product of two non-flat functions.
Thus, the formal Taylor series around $t_0$ of the left hand side
is identically zero, whereas the corresponding series
for the right hand side is non-zero.
\qed

Lemma \ref{lemma:D1} shows how to obtain $\omega_m(u, \cdot)$ directly from $u$,
without referring to $v(u,\cdot)$. More generally, 
we define the local $m$-argument $\hat\omega_m$ by the equation
\[ \hat\omega_m'(u,t) =
m - \frac{m^2 + f'(u(t))}{m} \sin^2(\hat\omega_m(u,t)),
\qquad \hat\omega_m(u,t_0) = \theta_0\]
where $\tan(\theta_0) = ma_0/b_0$.

\begin{prop}
\label{prop:D4}
Let $f$ be appropriate;
then $C_m \ne \emptyset$ if and only if the number
$- m^2$  belongs to the interior of the image of $f'$.
\end{prop}

Notice that the proposition makes no claim concerning
connectedness or any other topological property of $C_m$.

{\nobf Proof: }
Assume that $- m^2 \in \interior(f'(\RR))$:
for some $\epsilon > 0$ there exist $x_-, x_+ \in \RR$
with $f'(x_\pm) = - m^2 \mp \epsilon$. 
Consider two families of functions $u_{-,s}$, $u_{+,s}$
in $C^0_0([0,\pi])$ which are uniformly $C^0$ bounded
and $L^1$ converge to the constants $x_-, x_+$ when $s$ tends to 0.
From item (c) of Lemma \ref{lemma:D1}, 
for some sufficiently small $s_0$ we have
$\omega_m(u_{-,s_0},\pi) < -m^2$,
$\omega_m(u_{+,s_0},\pi) > -m^2$.
By continuity, some $u$ in the line segment from
$u_{-,s_0}$ to $u_{+,s_0}$ belongs to $C_m$.

To prove the converse, assume first $f'(0) = -m^2$:
from the definition of appropriateness
we must then have $f''(0) \ne 0$
and then clearly $- m^2 \in \interior(f'(\RR))$.
Assume instead $f'(0) \ne -m^2$, say $f'(0) > -m^2$,
and let $u \in C_m$.
For $t_-$ near 0 or $\pi$ we have $f'(u(t_-)) > -m^2$
and from Lemma \ref{lemma:D1} $\omega_m'(u,t_-) < m$
for such $t_- \ne 0, \pi$. 
We must then have $\omega_m'(u,t_+) > m$ for some $t_+ \in [0,\pi]$.
Thus $f'(u(t_+)) < - m^2$ and the result follows.
\qed

We present a rough sketch of the construction of the deformation
within $C_m$ of a family of functions $u$
to a final point $u_\ast$, as claimed in Proposition \ref{prop:D7}.
Actually, from Theorem \ref{theo:A},
it suffices that this deformation be continuous
in the $C^0$ norm instead of the $H^2$ norm.
A natural candidate for $u_\ast$ would be a constant function
equal to an arbitrary $x_m$ where $f'(x_m) = -m^2$:
in this case, $\omega_m(u,t) = mt$.
Unfortunately, this function does not satisfy the Dirichlet
boundary conditions: given the family of functions $u$ to be deformed,
we construct a fixed $u_\ast \in C_m$
which is constant equal to $x_m$ in a large interval $[a, \pi - a]$.

\begin{figure}[ht]
\begin{center}
\epsfig{height=6cm,file=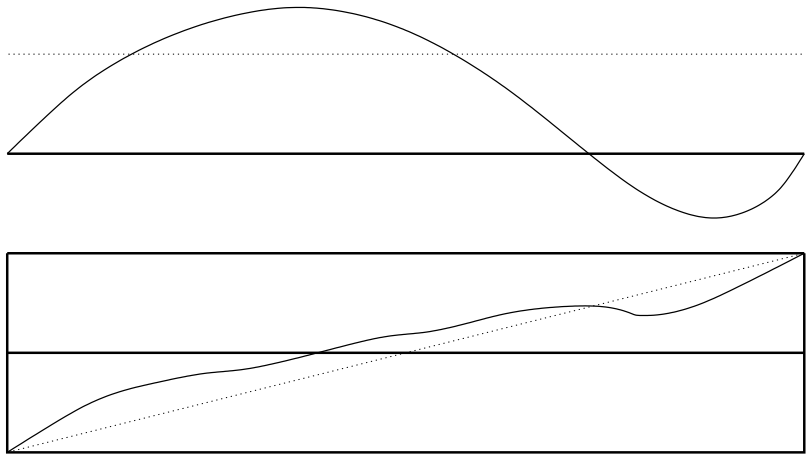}
\end{center}
\caption{Graph of original $u$ and $\omega_m$}
\label{fig:H0}
\end{figure}

\begin{figure}[ht]
\begin{center}
\epsfig{height=6cm,file=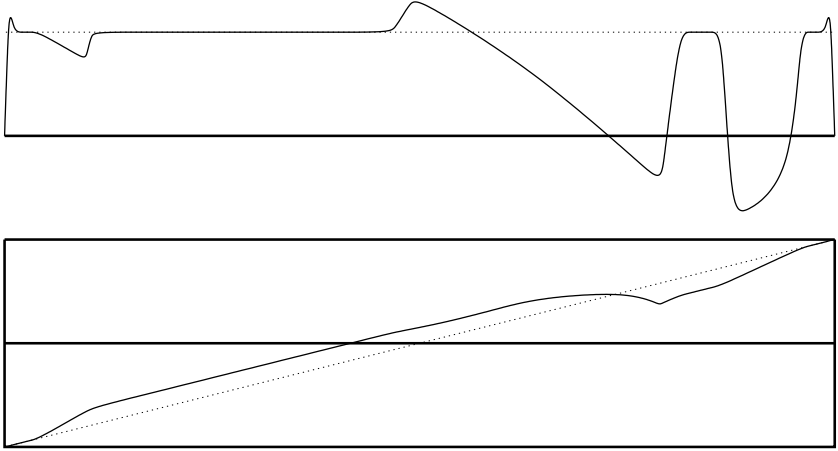}
\end{center}
\caption{Graph of $u$ and $\omega_m$ at some point during step 4}
\label{fig:H5}
\end{figure}

Consider now the $m$-argument of a given $u \in C_m$:
$\omega_m(u,\cdot)$ is a continuous function from $[0,\pi]$ to $[0, m\pi]$.
These two graphs are shown in Figure \ref{fig:H0};
in this example $m = 2$.
The graphs of the constant functions $u=x_m$
and $\omega_m(u,t)= tm$ are indicated by dotted lines.
The idea, to be formalized and implemented in step 4 of the proof,
is for the homotopy to squeeze the graph of $\omega_m$
between parallel walls advancing towards the line $y = mt$,
as shown in Figure \ref{fig:H5}.
A corresponding $u$ is obtained by changing its original value
in the region of the domain over which the graph of $\omega_m$
has been squeezed---there, the new value of $u$ is $x_m$.
Notice that, in principle, the value of $\omega_m(u,\pi) = m \pi$
for this new $u$ but such $u$ is discontinuous.
We must therefore make amends:
for a fixed tolerance $\tol$, the region where the graph of $\omega_m$
trespasses the wall by more than $\tol$ is taken to $x_m$
and in the region where the graph of $\omega_m$ lies strictly
between the walls, $u$ is unchanged.
Hence, there is an open region in the domain
where $u$ assumes rather arbitrary values
in order to preserve its continuity.
Steps 2 and 3 guarantee that this open region is uniformly small
(for $u$ in the deformed loop): in particular,
we have to fudge $u$ so that the graph of $\omega_m(u,\cdot)$
does not include segments parallel to $y = mt$.
This {\sl symmetry breaking} is achieved by replacing $u$
by appropriate polynomial approximations.
At the end of step 4, we have functions which are equal to $x_m$
in a substantial amount of the domain:
the straight line segments joining them to $u_\ast$
can be deformed onto $C_m$ thus accomplishing
the last step of the deformation.

\begin{lemma}
\label{lemma:D2}
Given $u \in C^0_D$, $\omega_m(u,\cdot)$ is a $C^1$ function with
$\omega'_m(u,t_0) = m$ and $\omega''_m(u,t_0) = 0$
whenever $\omega_m(u, t_0) = j\pi$, $j \in \ZZ$.
Furthermore, given $u_0 \in \RR$ with $f''(u_0) \ne 0$ and
a $C^1$ function $w: (t_0 - \epsilon, t_0 + \epsilon) \to \RR$
with either $w(t_0) \ne j\pi$ or simultaneously
$w(t_0) = j\pi$, $w'(t_0) = m$, $w''(t_0) = 0$,
there exists, for sufficiently small $\epsilon_1 < \epsilon$,
a unique continuous function
$u: (t_0 - \epsilon_1, t_0 + \epsilon_1) \to \RR$ with
\[ w'(t) =
m - \frac{m^2 + f'(u(t))}{m} \sin^2(w(t)),
\qquad u(t_0) = u_0.  \]
\end{lemma}

{\nobf Proof: }
The fact that $\omega_m(u,\cdot)$ is a $C^1$ function
follows directly from the differential equation in Lemma \ref{lemma:D1}.
Given $t_0$ with $\omega_m(u,t_0) = j\pi$, we have
\begin{align}
\omega''_m(u,t_0) &=
\lim_{t \to t_0} \frac{\omega'_m(u,t) - m}{t - t_0} \notag\\
&=
\lim_{t \to t_0} - \frac{m^2 + f'(u(t))}{m}\;
\frac{\sin^2(\omega_m(u,t))}{t - t_0}
\notag\\ &=
C \lim_{t \to t_0}
\frac{\sin^2(\omega_m(u,t)) - \sin^2(\omega_m(u,t_0))}{t - t_0} = 0. \notag
\end{align}
Also, the differential equation in $u$ and $w$ may be written as
\[ f'(u(t)) = - m^2 + m \frac{m - w'(t)}{\sin^2(w(t))}. \]
The right hand side is clearly continuous even when $w(t) = j\pi$
which allows for solving in $u$ when $f'$ is invertible.
\qed

Thus, for example, for concave or convex nonlinearities $f$,
changes in $\omega_m(u,\cdot)$ easily translate back to changes in $u$.
In particular, as mentioned in the introduction,
simpler proofs of Theorem \ref{theo:INTRO} are known
under these hypothesis.


A technical difficulty in the construction of $u_\ast$
(and in many points of the deformation process)
is making sure that $u_\ast$ (or a deformed $u$) belongs to $C_m$.
For this, we make use of the following Lemma.

Consider two smooth functions $u_0: [0, t_0] \to \RR$ and 
$u_1: [t_1, \pi] \to \RR$ satisfying $u_0(0) = 0$ and $u_1(\pi) = 0$, 
and consider local $m$-arguments
$\omega_m(u_0,t)$ and $\omega_m(u_1,t)$ 
starting respectively from $0$ at $t = 0$ and $m\pi$ at $t = \pi$.
Under the hypothesis of the lemma below, we may {\sl solder} 
these two chunks in a specific, smoothly defined way, to obtain
a smooth function $u:[0,\pi] \to \RR$ belonging to $C_m$ (i.e., so that
$\omega_m(u,\pi) = m\pi)$.

\begin{lemma}
\label{lemma:D3}
Let $0 < t_0 < t_1 < \pi$ be real numbers such that
$\sin(mt) \ne 0$ for all $t \in [t_0, t_1]$.
For sufficiently small $\epsilon > 0$ there exists a smooth
function 
\[ \Xi_{t_0, t_1}: (-\epsilon, \epsilon) \times (-\epsilon, \epsilon)
\times [t_0, t_1] \to \RR \]
with the following properties:
\begin{enumerate}[(a)]
\item{ $\Xi_{t_0, t_1}(h_0, h_1, t_0) = \Xi_{t_0, t_1}(h_0, h_1, t_1) = x_m$
for all $h_0, h_1 \in (-\epsilon, \epsilon)$ and flat at these points; }
\item{ if $\hat\omega_m(\Xi_{t_0, t_1}(h_0, h_1, \cdot), t)$ is defined
with initial condition
$\hat\omega_m(\Xi_{t_0, t_1}(h_0, h_1, \cdot), t_0) = m t_0 + h_0$
then $\hat\omega_m(\Xi_{t_0, t_1}(h_0, h_1, \cdot), t_1) = m t_1 + h_1$; }
\item{ $\Xi_{t_0,t_1}(h, h, t) = x_m$
for all $h \in (-\epsilon, \epsilon)$, $t \in [t_0, t_1]$.}
\end{enumerate}
\end{lemma}

{\nobf Proof: }
Take $\epsilon > 0$ such that 
$\sin(mt) \ne 0$ for all $t \in [t_0 - \epsilon/m, t_1 + \epsilon/m]$.
Let $\xi: [t_0, t_1] \to [0, 1]$ be a smooth bijection
with $\xi(t_0) = 0$, $\xi(t_1) = 1$ 
which is flat at these endpoints.
Let
\[ w_{h_0, h_1}(t) = mt + h_0 + (h_1 - h_0) \xi(t); \]
we are ready to apply Lemma \ref{lemma:D2} in order to define $\Xi_{t_0, t_1}$
with
\[ \hat\omega_m(\Xi_{t_0, t_1}(h_0, h_1, \cdot), t) = w_{h_0, h_1}(t). \]
\qed

The lemma below will be used in step 4
of the proof of Proposition \ref{prop:D7}.

\begin{lemma}
\label{lemma:D5}
Let $g: \Ss^k \times [0,\pi] \to \RR$ be a smooth function.
Suppose there exists $\epsilon > 0$ such that for any $\delta > 0$
there exists $y_\delta \in \RR$ and $\theta_\delta \in \Ss^k$
such that $\lebesgue(A_\delta) > \epsilon$, where
\[ A_\delta = \{ t \in [0,\pi] \;|\;
g(\theta_\delta, t) \in [y_\delta, y_\delta + \delta] \}. \]
Then there exists a point $(\theta_0, t_0) \in \Ss^k \times [0, \pi]$
for which the function $g(\theta_0, \cdot)$ is flat at $t_0$.
\end{lemma}

{\nobf Proof: }
Let $(y_0, \theta_0)$ be an accumulation point of the sequence
$(y_{a_n}, \theta_{a_n})$ where $\lim a_n = 0$;
without loss, we may suppose that $(y_0, \theta_0)$ is the limit
of this sequence.  Let 
\[ A_0 = \limsup_n A_{a_n} = \bigcap_n \bigcup_{i \ge n} A_i; \]
from
\[ \lebesgue(\limsup_n A_{a_n}) \ge \limsup_n \lebesgue(A_{a_n}) \]
(a corollary of Fatou's lemma)
we have $\lebesgue(A_0) \ge \epsilon$.
On the other hand, from continuity,
$g(\theta_0, t) = y_0$ for $t \in A_0$.
Thus, any non-isolated point $t_0 \in A_0$
is a flat point for $g(\theta_0, \cdot)$. \qed

{\nobf Proof of Proposition \ref{prop:D7}: }
Path connectedness of $C_m$
(i.e., the case $k = 0$, for which $\pi_k(C_m)$ has no group structure)
will be discussed simultaneously to the verification of the triviality
of the homotopy groups of $C_m$.
Keeping with the notation of the previous section, let
$\phi: X \to \RR$, $\phi(u) = \omega_m(u,\pi) - m\pi$,
$M_Y = C_m$ and $M_X = \phi^{-1}(\{0\}) \subseteq X$.
As shown in the previous section,
it suffices, given a $Y$-continuous (i.e., $H^2$-continuous)
$\gamma: \Ss^k \to M_Y$,
to construct a $X$-continuous (i.e., $C^0$-continuous)
$\Gamma: \BB^{k+1} \to M_X$
with $\Gamma|_{\Ss^k = \partial \BB^{k+1}} = \gamma$.

We set a more convenient notation:
let $A^{k+1}$ be the annulus $\Ss^k \times [0,\pi]$
and $U(0): A^{k+1} \to \RR$ be the initial map
defined by $U(0)(\theta,t) = (\gamma(\theta))(t)$.
We want to construct a continuous function $U: [0,5] \times A^{k+1} \to \RR$
with the following properties:
\begin{enumerate}[(i)]
\item{the values at $\{0\} \times A^{k+1}$ are given by the initial map $U(0)$;}
\item{for any $s \in [0,5]$ and $\theta \in \Ss^k$,
$U(s,\theta,\cdot) \in M_X \subseteq C^0_0([0,\pi])$, i.e.,
$U(s,\theta,0) = U(s,\theta,\pi) = 0$ and
$\omega_m(U(s,\theta,\cdot),\pi) = m\pi$;}
\item{for $s = 5$ the function $U$ is constant in $\theta$, i.e.,
for any $\theta_0, \theta_1 \in \Ss^k$ and $t \in [0,\pi]$,
$U(5,\theta_0,t) = U(5,\theta_1,t)$.}
\end{enumerate}
We then have
\[ \Gamma(r \theta)(t) = 
\begin{cases}
U(5, \theta, t), & \textrm{ if } r \le 1/2, \notag\\
U(10(1-r), \theta, t) & \textrm{ if } r \ge 1/2. \notag
\end{cases}
\]

\medskip

{\nobf Step 1 }
We search for convenient vectors $\varphi(\theta,\cdot)$ in $Y$
along which the derivative
\[ 
- \frac{m}{(v'(U(0,\theta,\cdot),\pi))^2}
\int_0^\pi f''(U(0,\theta,\sigma)) \varphi(\theta,\sigma)
(v(U(0,\theta,\cdot),\sigma))^2 d\sigma, \]
of $\omega_m(U(0,\theta,\cdot),\pi)$ is positive
(see Lemma \ref{lemma:D1}, (d)).
Let $\beta_\delta: [0,\pi] \to [0,1]$ be a smooth bump
equal to zero in $[0, \delta) \cup (\pi - \delta, \pi]$
and equal to one in $(2\delta, \pi - 2\delta)$;
this family of bumps may be constructed to be $L^1$-continuous in $\delta$.
Clearly, the choice
$\varphi(\theta,t) = - \beta_\delta(t) f''(U(0,\theta,t))$
yields a positive derivative for $\delta = 0$
and therefore, by continuity, also for some positive $\delta_0$
(independent of $\theta$).
Furthermore, the smoothness in $u$ of $\omega_m(u,\pi)$
guarantees that there exists $\epsilon > 0$
such that
\[ \frac{\partial}{\partial \tau}
\omega_m(U(0,\theta,\cdot) + \tau \varphi(\theta, \cdot), \pi) > \epsilon \]
at any point $(\theta, \tau)$ with $|\tau| < \epsilon$.

We now want to deform each $u$ so that the modified $u$'s are constant
equal to $x_m$ on small intervals near $t = 0$ and $t = \pi$.
More precisely, we will define small positive real numbers
$0 < \delta_2 \ll \delta_1$ and construct $U$ for $0 \le s \le 1$ so that
\begin{enumerate}[(i)]
\item { for $\delta_2 \le t \le \delta_1$ and
$\pi - \delta_1 \le t \le \pi - \delta_2$ we have
$U(1,\theta,t) = x_m$; }
\item { for $t < \delta_2$ and $t > \pi - \delta_2$,
$U(1,\theta,t)$ does not depend on $\theta$
and satisfies $|U(1,\theta,t)| \le |x_m|$; }
\item { for any $s \in [0,1]$, and any $\theta \in \Ss^k$
we have $\omega_m(U(s,\theta,\cdot),\pi) = m\pi$,
i.e., $U(s,\theta,\cdot) \in C_m$. }
\end{enumerate}
In particular, from {(i)},
in the interval $\delta_2 \le t \le \delta_1$
(or $\pi - \delta_1 \le t \le \pi - \delta_2$)
the graph of $\omega_m(U(1,\theta,\cdot),t)$ is a straight segment
of slope $m$.
We claim that for any sufficiently small $\delta_1$ and $\delta_2$
with $\delta_2 < \delta_1$ this construction can be accomplished.
Indeed, for sufficiently small $\delta_1 < \delta_0/2$
we can uniquely define
$\xi: [0,1] \times \Ss^k \times [0,\pi] \to (-\epsilon, \epsilon)$
such that
\[ U(s,\theta,t) = (1 - s + s\beta_{\delta_1}(t)) U(0,\theta,t) +
s (\beta_{\delta_2/2} - \beta_{\delta_1})(t) x_m +
\xi(s,\theta,t) \varphi(\theta,t) \]
satisfies $U(s,\theta,\cdot) \in C_m$ for all $s \in [0,1]$
and $\theta \in \Ss^k$
(here $\delta_0$, $\epsilon$, $\varphi$ and $\beta$ are
as defined in the previous paragraph).
The above formula for $U$ admits the following geometric interpretation.
The first two terms of the right hand side parametrize in $s$
a straight line segment from the original $u$
to a modified function satisfying items {(i)} and {(ii)}.
The third term takes care of item {(iii)} provided 
the $L^1$ distance between $u$ and the modified function
is so small that their $m$-arguments $\omega_m$ at $\pi$
differ from less than $\epsilon$.

\smallskip

\begin{figure}[ht]
\begin{center}
\epsfig{height=6cm,file=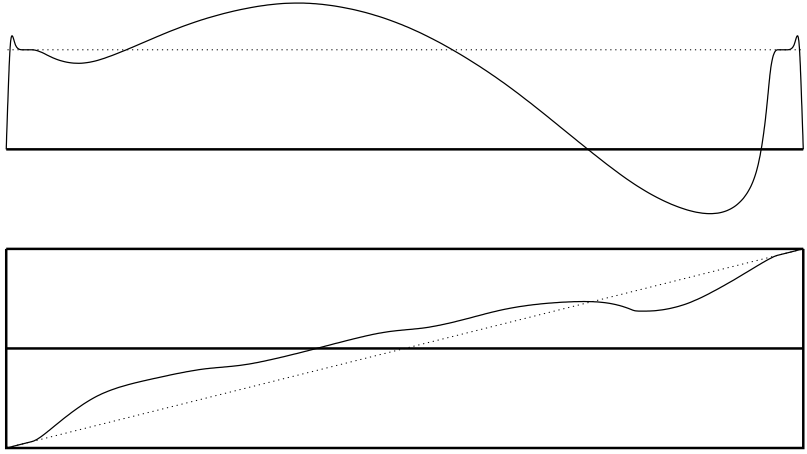}
\end{center}
\caption{Graphs of $u$ and $\omega_m$ at the end of step 2}
\label{fig:H1}
\end{figure}

{\nobf Step 2}
Let $\eta > 0$ be a small number to be specified later.
We now define $U$ for $1 \le s \le 2$ so that
there exist positive numbers $\delta'_2 < \delta'_1 < \delta'_0 < \pi/m$
such that
\begin{enumerate}[(i)]
\item { for $t \in [\delta'_2, \delta'_1] \cup
[\pi - \delta'_1, \pi - \delta'_2]$ we have
$U(2,\theta,t) = x_m$ and $\omega_m(U(2,\theta,\cdot),t) = mt$; }
\item { for any $\theta$, $U(2,\theta,\pi - \delta'_0) = x_m$ and
$\omega_m(U(2,\theta,\cdot), \pi - \delta'_0)) = m(\pi - \delta'_0) + \eta$; }
\item { for $t \in [0, \delta'_1] \cup [\pi - \delta'_0, \pi]$,
$U(2,\theta,t)$ does not depend on $\theta$; }
\item { for any $s \in [1, 2]$, and any $\theta \in \Ss^k$,
we have $U(s,\theta,\cdot) \in C_m$. }
\end{enumerate}
To do this, assume that $\delta_1 < \pi/m$ and $\delta_2 < \delta_1/4$
and let $\delta'_2 = \delta_1/2$, $\delta'_1 = 3\delta_1/4$,
$\delta'_0 = 7\delta_1/8$.
The $L^1$-continuity of the $m$-argument (Lemma \ref{lemma:D1})
entails that both
$h_- = m \delta_1/4 - \omega_m(U(1,\theta,\cdot), \delta_1/4)$ and
$h_+ = m (\pi - \delta_1/4) - \omega_m(U(1,\theta,\cdot), \pi - \delta_1/4)$
can be taken to be arbitrarily small by choosing $\delta_2$ small.
Now solder the six chunks: for $s \in [1, 2]$ set
\[ U(s,\theta,t) = \begin{cases}
\Xi_{\delta_1/4, \delta'_2}(h_-, h_-(s), t),
& \textrm{if } \delta_1/4 < t < \delta'_2 \notag\\
\Xi_{\delta'_1, \delta_1}(h_-(s), h_-, t),
& \textrm{if } \delta'_1 < t < \delta_1 \notag\\
\Xi_{\pi - \delta_1, \pi - \delta'_0}(h_+, h^0_+(s), t),
& \textrm{if } \pi - \delta_1 < t < \pi - \delta'_0 \notag\\
\Xi_{\pi - \delta'_0, \pi - \delta'_1}(h^0_+, h^1_+(s), t),
& \textrm{if } \pi - \delta'_0 < t < \pi - \delta'_1 \notag\\
\Xi_{\pi - \delta'_2, \pi - \delta_1/4}(h^1_+(s), h_+, t),
& \textrm{if } \pi - \delta'_2 < t < \pi - \delta_1/4 \notag\\
U(1/4, \theta, t), & \textrm{otherwise;} \notag
\end{cases} \]
where $h_-(s) =  (2 - s) h_-$, $h^0_+(s) = \eta +  (2 - s) (h_+ - \eta)$
and $h^1_+(s) =  (2 - s) h_+$.
The number $\eta$ is chosen to be any positive number
sufficiently small to permit soldering.
See Figure \ref{fig:H1} for a sketch of the graph of $U(2,\theta,\cdot)$ 
and $\omega_m(U(2,\theta,\cdot),t)$.

We are now ready to describe $u_\ast$:
\[ u_\ast(t) = \begin{cases}
U(2, \theta, t), & \textrm{ for }
t \in [0, \delta'_1] \cup [\pi - \delta'_1, \pi], \notag\\
x_m, & \textrm{ for } t \in [\delta'_2, \pi - \delta'_2]. \notag
\end{cases} \]
Notice that the value of $\theta$ in the definition above
is immaterial, from property {(iii)}.
Moreover, the homotopy from now on will not change the values
of the functions $U(s, \theta, \cdot)$ near the endpoints:
we shall have $U(s, \theta, t) = u_\ast(t)$
for all $s \ge 2$, $t \in [0, \delta'_2] \cup [\pi - \delta'_2, \pi]$
and $\theta \in \Ss^k$.
The construction of $u_\ast$ is dependent on the original 
map $\gamma: \Ss^k \to M_Y$:
still, path connectivity of $C_m$ follows from the same homotopy
being described when applied to an original map $\gamma$, with $k=0$.

\smallskip

{\nobf Step 3}
Continuing with the preparations, we change $u$'s
in the central interval $[\delta'_1, \pi - \delta'_0]$
so that there they become polynomials in the variables $\theta$ and $t$.
More precisely, for some polynomial function $P: \RR^{k+2} \to \RR$
(with properties to be specified below)
and for $s \in [2, 3]$, set
\[
U(s, \theta, t) =
\begin{cases}
(s - 2) U(2, \theta, t) + (3 - s) P(\theta, t), 
& \textrm{for } t \in [\delta'_1, \pi - \delta'_0], \notag\\
\Xi_{\pi - \delta'_1, \pi - \delta'_2}(h(s,\theta),0,t),
& \textrm{for } t \in [\pi - \delta'_1, \pi - \delta'_2], \notag\\
u_\ast(t), & \textrm{otherwise,}
\end{cases}
\]
where
$h(s, \theta) =
\omega_m(U(s,\theta, \cdot), \pi - \delta'_1) - m(\pi - \delta'_1)$.
In order for this function $U$ to be continuous we must choose $P$
with $P(\theta, \delta'_1) = P(\theta, \pi - \delta'_0) = x_m$.
Lemma \ref{lemma:D3} applies if we require $h(s, \theta) \in (\eta/2, 3\eta/2)$
for all $s$ and $\theta$:
this is accomplished from Lemma \ref{lemma:D1} by choosing
$P(\theta, t)$ uniformly
close to $U(2, \theta, t)$ in the interval $[\delta'_1, \pi - \delta'_0]$
(as we may, from the Stone-Weierstrass theorem).

\smallskip

{\nobf Step 4}
Let
\[ A = [3,4] \times \Ss^k \times [\delta'_1, \pi - \delta'_1] \] 
and consider the partition
\begin{align}
A_I &= \{ (s, \theta, t) \in A \;|\;
- (4-s)m\pi \le \omega_m(U(s, \theta, \cdot), t) - mt \le (4-s)m\pi \}, \notag\\
A_S &= \{ (s, \theta, t) \in A \;|\;
- (4-s)m\pi - \tol \ge \omega_m(U(s, \theta, \cdot), t) - mt \notag\\
&\qquad\qquad\qquad \textrm{or}\
\omega_m(U(s, \theta, \cdot), t) - mt \ge (4-s)m\pi + \tol \}, \notag\\
A_T &= A - (A_I \cup A_S) \notag
\end{align}
(here the letters $I, S$ and $T$ stand for {\sl invariant},
{\sl squeezed} and {\sl Tietze}).
For $(s, \theta, t) \in A$, set
\[ U(s, \theta, t) = \begin{cases}
U(3, \theta, t), & \textrm{ for } (s, \theta, t) \in A_I, \notag\\
x_m, & \textrm{ for } (s, \theta, t) \in A_S, \notag \end{cases} \]
and now apply the Tietze extension theorem to define $U$ on $A_T$
so that $U$ is continuous in $A$.
Set $U(s, \theta, t) = u_\ast(t)$ for
$t \in [0, \delta'_1] \cup [\pi - \delta'_2, \pi]$
and apply solder in the interval $[\pi - \delta'_1, \pi - \delta'_2]$,
i.e., set
\[ U(s, \theta, t) = \Xi_{\pi - \delta'_1, \pi - \delta'_2}%
(\omega_m(U(s, \theta, \cdot), \pi - \delta'_1), 0, t). \]
We are left to show that under this construction,
for sufficiently small $\tol$, solder is applicable, i.e.,
$\omega_m(U(s, \theta, \cdot), \pi - \delta'_1)$
can be taken as small as desired.

We first prove that given $\epsilon > 0$ there exists $\tol > 0$
such that for any $s_0$ and $\theta_0$ the Lebesgue measure of
\[ A_{T,s_0,\theta_0} = \{t \in [0,\pi] \;|\; (s_0, \theta_0, t) \in A_T\} \]
is smaller than $\epsilon$.
First consider
$A_{T,s_0,\theta_0} \cap [\pi - \delta'_0, \pi - \delta'_1]$: in this interval,
$\omega_m(U(3, \theta_0, \cdot), t) - mt$
is strictly decreasing and does not depend on $\theta_0$.
Thus, for sufficiently small $\tol$ we may assume
$\lebesgue(A_{T,s_0,\theta_0} \cap [\pi - \delta'_0, \pi - \delta'_1])
< \epsilon/2$.
Next consider $A_{T,s_0,\theta_0} \cap [\delta'_1, \pi - \delta'_0]$.
Suppose by contradiction that for a fixed $\epsilon > 0$
and for all $\tol > 0$ we have
\[ \lebesgue(A_{T,s_0,\theta_0} \cap [\delta'_1, \pi - \delta'_0])
> \epsilon/2. \]
From Lemma \ref{lemma:D5},
$\omega'_m(U(3,\theta_0,\cdot),t)$ is flat in $t = t_0$
for some $\theta_0 \in \Ss^k$.
Now, making use of the last item of Lemma \ref{lemma:D1},
$U(3, \theta_0, t)$ is flat in $t = t_0$,
contradicting its polynomiality.

\goodbreak
\smallskip

{\nobf Step 5}

At this point of the construction the loop $U(4,\theta,t)$ is
$L^1$ close to the constant loop $u_\ast$.
Indeed, consider $U(3,\theta,t)$:
this function of $\theta$ and $t$ is $C^0$-bounded by some constant $C$.
Given $\epsilon_1 > 0$ we can choose $\tol$ in step 4
in such a way that, for all $\theta$,
\begin{enumerate}
\item{ $|U(4,\theta,\cdot) - u_\ast|_{C^0} < 2 C + 1$, }
\item{ $|U(4,\theta,\cdot) - u_\ast|_{L^1} < \epsilon_1$. }
\end{enumerate}
The first follows from compactness of the set of parameters
for the function $\Xi$ employed in the soldering in step 4.
The second follows from choosing $\tol$ such that $A_T$ has measure
less than $\epsilon_1/(4C + 2)$ (as discussed in step 4).

We may now take a family in $\theta$ of straight lines 
joining $U(4,\theta,t)$ to $U(5,\theta,t) = u_\ast$.
More formally, let
\[
U(s, \theta, t) =
\begin{cases}
\Xi_{\pi - \delta'_1, \pi - \delta'_2}(h(s,\theta),0,t),
& \textrm{for } t \in [\pi - \delta'_1, \pi - \delta'_2], \notag\\
(5 - s) U(4,\theta,t) + (s - 4) u_\ast(t),  & \textrm{otherwise,}
\end{cases}
\]
where, again,
$h(s, \theta) =
\omega_m(U(s,\theta, \cdot), \pi - \delta'_1) - m(\pi - \delta'_1)$.
The fact that $U(s,\theta,\pi - \delta'_1) - m(\pi - \delta'_1)$
is appropriately small follows from Lemma \ref{lemma:D1}, item (c).
\qed


\bigskip\bigskip\bigbreak

{

\parindent=0pt
\parskip=0pt
\obeylines

Dan Burghelea, Department of Mathematics, Ohio State University,
231 West 18th Ave, Columbus, OH 43210-1174, USA

\smallskip

burghele@math.ohio-state.edu

\smallskip

Nicolau C. Saldanha and Carlos Tomei, Departamento de Matem\'atica, PUC-Rio
R. Marqu\^es de S. Vicente 225, Rio de Janeiro, RJ 22453-900, Brazil

\smallskip

nicolau@mat.puc-rio.br; http://www.mat.puc-rio.br/$\sim$nicolau/
tomei@mat.puc-rio.br

}

\end{document}